\documentclass[english,10pt,a4paper]{article}
\usepackage{blindtext}
\usepackage[T1]{fontenc}
\usepackage[utf8]{inputenc}
\usepackage{graphicx}
\usepackage{babel}
\usepackage{amsmath,amsthm,amssymb}
\usepackage{bbold}
\usepackage[cyr]{aeguill}
\usepackage{a4}
\usepackage{fancyhdr}
\usepackage{lastpage}
\usepackage{hyperref}
\usepackage{csquotes}
\usepackage{enumitem}
\usepackage{pifont}
\usepackage{comment}
\usepackage{mathtools}
\usepackage{authblk}
\usepackage{stmaryrd}
\usepackage{bbm}
\usepackage{hyperref}
\usepackage{multirow}
\usepackage{tabularx}
\usepackage[Bjornstrup]{fncychap}
\usepackage{resmes}
\usepackage{mathrsfs}
\usepackage{minted}
\usepackage[nottoc,numbib]{tocbibind}
\usepackage[all]{xy}
\usepackage{tikz}
\usepackage{tikz-cd}
\usepackage{authblk}
\usepackage[multiple]{footmisc}

\usetikzlibrary{positioning}

\usepackage[backend=biber,style=numeric]{biblatex}
\renewbibmacro{in:}{\ifentrytype{article}
    {}
    {\bibstring{in}\printunit{\intitlepunct}}
}

\newtheorem{lemma}{Lemma}
\newtheorem{definition}{Definition}
\newtheorem{proposition}{Proposition}
\newtheorem{theorem}{Theorem}
\newtheorem{remark}{Remark}
\newtheorem*{remark*}{Remark}
\newtheorem{corollary}{Corollary}

\makeatletter
\renewcommand{\@fnsymbol}[1]{\@arabic{#1}}
\makeatother

\title{Tensor Products with Verma Module and Restriction to Parabolic Subalgebra}

\author{
Antoine Merceron\thanks{CMLS, École polytechnique, France and Mathematisches Institut, Albert-Ludwigs-Universität Freiburg, Germany.}
$^{,}$\thanks{I thank Wolfgang Soergel for his guidance and feedback during the preparation of this article.}
}

\date{September 17, 2025}

\begin{document}
\maketitle

\begin{abstract}
\indent Using the tensor identity, we obtain decomposition results for the tensor product of a generalized Verma module with a module $M$ in the category $\mathcal{O}^{\mathfrak{p}}$, based on the decomposition of the restriction of $M$ to the parabolic subalgebra $\mathfrak{p}$. We show that this restriction admits an essentially unique decomposition into indecomposable $\mathfrak{p}$-modules and identify two particular types of direct summands: finite-dimensional quotients and tilting submodules. Finally, we give the complete $\mathfrak{b}$-decomposition of all indecomposable $M \in \mathcal{O}$ in $\mathfrak{sl}_2$, and of all Verma modules in the block of a dominant integral weight in $\mathfrak{sl}_3$, from which we derive explicit computations.
\end{abstract}

\section{Introduction}

Let $\mathfrak{g}$ be a finite-dimensional complex semisimple Lie algebra with a Cartan decomposition $\mathfrak{g} = \mathfrak{n}^- \oplus \mathfrak{h} \oplus \mathfrak{n}$ where $\mathfrak{b} = \mathfrak{h} \oplus \mathfrak{n}$ is a Borel subalgebra. We also fix a parabolic subalgebra $\mathfrak{p} = \mathfrak{p}_I \supset \mathfrak{b}$. Following the definitions and notations of Humphreys \cite{Humphreys2008}, we consider the BGG category $\mathcal{O}$ as well as the parabolic subcategory $\mathcal{O}^{\mathfrak{p}}$. \\

For $\lambda \in \mathfrak{h}^*$, the Verma module of highest weight $\lambda$ is defined as $M(\lambda) = \operatorname{Ind}_{\mathfrak{b}}^{\mathfrak{g}}\bigl(\mathbb{C}_\lambda\bigr) \;\in\; \mathcal{O}$. More generally, for $\lambda \in \Lambda_I^+$, we define the generalized Verma module of highest weight $\lambda$ as $M_{\mathfrak{p}}(\lambda) = \operatorname{Ind}_{\mathfrak{p}}^{\mathfrak{g}}\bigl(L_I(\lambda)\bigr) \;\in\; \mathcal{O}^{\mathfrak{p}}$. All of our results will be stated in the parabolic setting, though the reader may restrict attention to the special case $\mathfrak{p} = \mathfrak{b}$, for which $I = \emptyset, \quad \mathcal{O}^{\mathfrak{p}} = \mathcal{O}, \quad \Lambda_I^+ = \mathfrak{h}^*, \quad M_{\mathfrak{p}}(\lambda) = M(\lambda)$. \\

Our focus will be the decomposition of tensor products $M_{\mathfrak{p}}(\lambda) \otimes M$, with $M \in \mathcal{O}^{\mathfrak{p}}$. We begin in Section \ref{Section2} by recalling the tensor identity and observing that any direct summand of $\operatorname{Res}_{\mathfrak{p}}^{\mathfrak{g}} M$ gives rise to a direct summand of $M_{\mathfrak{p}}(\lambda) \otimes M$. In Section \ref{Section3}, we show that $\operatorname{Res}_{\mathfrak{p}}^{\mathfrak{g}} M$ admits an essentially unique decomposition into indecomposable $\mathfrak{p}$-modules. Using the adjunction $(\operatorname{Ind}_{\mathfrak{p}}^{\mathfrak{g}}, \operatorname{Res}_{\mathfrak{p}}^{\mathfrak{g}})$, we show that tilting submodules and finite-dimensional quotients are always direct summands of $\operatorname{Res}_{\mathfrak{p}}^{\mathfrak{g}} M$. Finally, in Section \ref{Section4}, we provide explicit decompositions of certain modules. In the case of $\mathfrak{sl}_2$, we decompose all indecomposable modules of category $\mathcal{O}$ and show that this allows us to deduce the decomposition of all tensor products $M(\lambda)  \otimes M$ with $M \in \mathcal{O}$. In the case of $\mathfrak{sl}_3$, we give the decomposition of $M(w\cdot\lambda)$ into indecomposable $\mathfrak{b}$-modules for $w \in W$ and $\lambda \in \Lambda^+$ a dominant integral weight. From this, we deduce the decomposition into indecomposables of $\bigl(M \otimes M'\bigr)^{\chi_0}$ for $M$ and $M'$ Verma modules in the principal block.

\section{Preliminaries}
\label{Section2}

We begin by recalling some elementary results on the tensor product of two modules in category $\mathcal{O}$, as presented in K{\aa}hrstr{\"o}m \cite{Kahrstrom2010}:

\begin{proposition}
\label{PropPreliStructure}
For $M, N \in \mathcal{O}^{\mathfrak{p}}$, the tensor product $M \otimes N$ satisfies the following properties:
\begin{enumerate}
\item $M \otimes N$ is a weight module. The dimensions of its weight spaces are finite and given by
$$
\dim (M \otimes N)_\lambda = \sum_{\lambda = \mu + \nu} \dim M_{\mu}\,\dim N_{\nu}
$$
\item $M \otimes N$ is locally $Z(\mathfrak{g})$-finite, and decomposes as
$$
M \otimes N = \bigoplus_{\chi} (M \otimes N)^{\chi}, \quad (M \otimes N)^{\chi} \in \mathcal{O}^{\mathfrak{p}}
$$
where $\chi$ runs over every central character $\chi : Z(\mathfrak{g}) \to \mathbb{C}$. In particular, $M \otimes N$ decomposes as a direct sum of indecomposable objects in $\mathcal{O}^{\mathfrak{p}}$, uniquely determined up to isomorphism and permutation of the summands. 
\item $(M \otimes N)^\vee \cong M^\vee \otimes N^\vee$
\end{enumerate}
\end{proposition}

\noindent
We also recall the tensor identity:

\begin{proposition}[Tensor identity]
For any $\mathfrak{g}$-module $M$ and $\lambda \in \Lambda_I^+$:
$$
M_{\mathfrak{p}}(\lambda) \otimes M \;\cong\; \operatorname{Ind}_{\mathfrak{p}}^{\mathfrak{g}}\bigl(L_I(\lambda) \otimes \operatorname{Res}_{\mathfrak{p}}^{\mathfrak{g}} M \bigr)
$$
In particular,
$$
M_{\mathfrak{p}}(0) \otimes M \;\cong\; \operatorname{Ind}_{\mathfrak{p}}^{\mathfrak{g}} \operatorname{Res}_{\mathfrak{p}}^{\mathfrak{g}} M
$$
\end{proposition}

A module $M \in \mathcal{O}^{\mathfrak{p}}$ is said to admit a generalized Verma filtration if it has a (finite) filtration whose successive quotients are isomorphic to generalized Verma modules $M_{\mathfrak{p}}(\mu)$. Moreover, $M \in \mathcal{O}^{\mathfrak{p}}$ is called a $\mathfrak{p}$-tilting module if both $M$ and $M^\vee$ admit generalized Verma filtration. We recall the three following facts (see \cite[§1.5-1.7]{Andersen1995}):
\begin{itemize}
\item A generalized Verma module $M_{\mathfrak{p}}(\mu)$ is $\mathfrak{p}$-tilting if and only if $M_{\mathfrak{p}}(\mu) = L(\mu)$.
\item If $M$ and $N^\vee$ admit a generalized Verma filtration then $\operatorname{Ext}_{\mathcal{O}^{\mathfrak{p}}}(M,N) = 0$.
\item The indecomposable $\mathfrak{p}$-tiltings are, up to isomorphism, of the form $D_{\mathfrak{p}}(\lambda)$ with $\lambda \in \Lambda_I^+$, where $D_{\mathfrak{p}}(\lambda)$ denotes the unique indecomposable $\mathfrak{p}$-tilting whose weights satisfy $\mu \leq \lambda$.
\end{itemize} 

\noindent
An important consequence of the tensor identity is the following result \cite[Prop.~2.6]{Kahrstrom2010}:

\begin{proposition}
\label{PropPreliFiltration}
For every $M, N \in \mathcal{O}^{\mathfrak{p}}$ and central character $\chi$:
\begin{itemize}
\item If $M$ admits a generalized Verma filtration, then $\bigl(M_{\mathfrak{p}}(\lambda) \otimes N\bigr)^\chi \in \mathcal{O}^{\mathfrak{p}}$ also admits a generalized Verma filtration.
\item If $M$ is a $\mathfrak{p}$-tilting, then $\bigl(M \otimes N\bigr)^\chi \in \mathcal{O}^{\mathfrak{p}}$ is also a $\mathfrak{p}$-tilting.
\item If $M$ and $N^\vee$ admit a generalized Verma filtration, then $\bigl(M \otimes N\bigr)^\chi \in \mathcal{O}^{\mathfrak{p}}$ is a $\mathfrak{p}$-tilting. 
\end{itemize}
\end{proposition}
In particular, by Propositions \ref{PropPreliStructure} and \ref{PropPreliFiltration}, the tensor product $M_{\mathfrak{p}}(\lambda) \otimes M$, $M \in \mathcal{O}^{\mathfrak{p}}$, admits an essentially unique decomposition into indecomposables with generalized Verma filtration in $\mathcal{O}^{\mathfrak{p}}$. \\

It follows from the tensor identity that every direct summand of $\operatorname{Res}_{\mathfrak{p}}^{\mathfrak{g}} M$ gives rise to a direct summand of $M_{\mathfrak{p}}(\lambda) \otimes M$. Thus, the problem of decomposing this tensor product is related to the problem of decomposing $M$ as a $\mathfrak{p}$-module. This strategy will be employed in the following sections to find direct summands of $M_{\mathfrak{p}}(\lambda) \otimes M$. \\

It should be emphasized, however, that this method does not yield the full indecomposable decomposition of $M_{\mathfrak{p}}(\lambda) \otimes M$. For example, if $M = L(\lambda)$ is simple, then it is indecomposable as a $\mathfrak{p}$-module. Nevertheless, the tensor product $M_{\mathfrak{p}}(\lambda) \otimes L(\lambda)$ typically decomposes into infinitely many indecomposable modules in category $\mathcal{O}$.

\section{$\mathfrak{p}$-module decomposition}
\label{Section3}

In this section, we present some general results on the restriction to $\mathfrak{p}$ of a module $M \in \mathcal{O}^\mathfrak{p}$.

\subsection{Decomposition into indecomposables}

If $\mathfrak{p} \neq \mathfrak{g}$, then for $M \in \mathcal{O}$, the restriction $\operatorname{Res}_{\mathfrak{p}}^{\mathfrak{g}} M$ is typically not noetherian, so Krull–Schmidt theorem does not apply a priori. Nonetheless, we will show that it still holds in this context. 

\begin{definition}
Let $M$ a $\mathfrak{g}$-module, for $\lambda \in \mathfrak{h}^*$, we denote by $M^+_\lambda$ the vector space of vectors $v^+ \in M_\lambda$ such that $\mathfrak{n} \cdot v^+ = 0$. Any nonzero vector $v^+$ of this form is called a maximal vector of weight $\lambda$. We denote by $\Pi^{+}(M) \subset \Pi(M)$ the set of weights of the maximal vectors of $M$.
\end{definition}

Let $\mathcal{M}_{\mathfrak{p}}$ be the additive full subcategory of the category of $\mathfrak{p}$-modules, whose objects $M$ satisfies the following four properties:
\begin{enumerate}
\item $M$ is $\mathfrak{h}$-semisimple
\item $M$ is locally $\mathfrak{n}$-finite
\item $\operatorname{dim} M_\lambda < \infty$ for all $\lambda \in \mathfrak{h}^*$
\item $\Pi^+(M)$ is finite
\end{enumerate}
In particular, $\mathcal{M}_{\mathfrak{p}}$ contains the restriction to $\mathfrak{p}$ of any module in category $\mathcal{O}$. \\

We will show that $\mathcal{M}_{\mathfrak{p}}$ is a Krull-Schmidt category, i.e. every object in $\mathcal{M}_{\mathfrak{p}}$ decomposes as a finite direct sum of indecomposables modules with local endomorphism rings (see \cite{Krause2015}). As a result, this decomposition is unique up to isomorphism and permutation of the summands.

\begin{lemma}
\label{LemDecomposition}
Let $(M_i)_{i \in I}$ be a family of $\mathfrak{p}$-submodules of $M \in \mathcal{M}_{\mathfrak{p}}$. The following two conditions are equivalent:
\begin{enumerate}
\item $\sum_{i \in I} (M_i)^+_\lambda$ is direct for all $\lambda \in \mathfrak{h}^*$,
\item $\sum_{i \in I} M_i$ is direct
\end{enumerate}
Moreover, in this case we have $\left(\bigoplus M_i\right)^+_\lambda = \bigoplus_{i \in I} (M_i)^+_\lambda \quad \text{for all } \lambda \in \mathfrak{h}^*$.
\end{lemma}

\begin{proof}
\mbox{}
\begin{enumerate}
\item Suppose $\sum_{i \in I} (M_i)^+_\lambda$ direct for all $\lambda \in \mathfrak{h}^*$. Assume, by contradiction, that for some $\lambda\in \mathfrak{h}^*$ there exists a nontrivial linear relation among elements of the weight space $(M_i)_\lambda$:
$$
\sum_{j \in J} v_j = 0
$$
where $J \subset I$ is finite and $v_j \in (M_j)_\lambda$ are nonzero. Choose $\mu \in \mathfrak{h}^*$ and $x \in U(\mathfrak{b})$, $i \in J$, such that $x \cdot v_i$ is a maximal vector of weight $\mu$. By finiteness, we may assume that $\mu$ is maximal with respect to $\leq$. Then for $x$ and $i$ as above:
$$
x \cdot v_i = -\sum_{j \neq i} x \cdot v_j
$$
For $j \neq i$, if $x \cdot v_j \neq 0$, there exists $x' \in U(\mathfrak{b})$ such that $x' x \cdot v_j$ is a maximal vector of weight $\nu \geq \mu$. By maximality of $\mu$, $\nu = \mu$, so $x \cdot v_j \in (M_j)^+_\mu$. This contradicts the existence of a nontrivial relation.

\item Conversely, if $\sum_{i \in I} M_i$ is direct, then $\sum_{i \in I} (M_i)^+_\lambda$ is direct for all $\lambda \in \mathfrak{h}^*$.  Moreover, for $v \in \left(\bigoplus M_i\right)^+_\lambda$, write $v = \sum v_i$ with $v_i \in (M_i)_\lambda$. For all $x \in \mathfrak{n}$:
$$
0 = x \cdot v = \sum x \cdot v_i \implies x \cdot v_i = 0
$$
hence $v_i \in (M_i)^+_\lambda$ for all $i \in I$. \qedhere
\end{enumerate}
\end{proof}
\noindent
Let $M, N \in \mathcal{M}_{\mathfrak{p}}$. Recall that if $f \in \operatorname{Hom}_{\mathfrak{b}}(M,N)$ and $\lambda \in \mathfrak{h}^*$, then $f(M_\lambda) \subset N_\lambda$ so that $f$ induces a linear map $f_\lambda \in \operatorname{Hom}_{\mathbb{C}}(M_\lambda,N_\lambda)$.

\begin{lemma}
\label{LemHom}
The following natural map is injective:
$$
\operatorname{Hom}_{\mathfrak{b}}(M,N) \hookrightarrow \prod_{\lambda \in \Pi^+(N)} \operatorname{Hom}_{\mathbb{C}}(M_\lambda,N_\lambda)
$$
\end{lemma}
\begin{proof}
Let $f \in \operatorname{Hom}_{\mathfrak{b}}(M,N)$ be nonzero, and pick $v \in M$ such that $f(v) \neq 0$. There exists $x \in U(\mathfrak{b})$ such that $x \cdot f(v)$ is a maximal vector in $N$. Then $x \cdot v \in M$ is mapped by $f$ to this maximal vector. Denote by $\lambda$ the weight of this maximal vector, the induced map $f_\lambda \in \operatorname{Hom}_{\mathbb{C}}(M_\lambda, N_\lambda)$ is nonzero.
\end{proof}

\begin{lemma}[Fitting]
\label{LemFitting}
Let $M \in \mathcal{M}_{\mathfrak{p}}$ and $f \in \operatorname{End}_{\mathfrak{b}}(M)$. Then, for some $n \geq 1$, we have:
$$
\operatorname{Ker}(f^n)\oplus \operatorname{Im}(f^n) = M
$$
\end{lemma}
\begin{proof}For each $\lambda \in \Pi^+(M)$, the weight space $M_\lambda$ is finite-dimensional, so the linear map $f_\lambda$ has a nonzero annihilating polynomial. Multiplying these polynomials over all $\lambda \in \Pi^+(M)$, we obtain a single polynomial that annihilates all $f_\lambda$ simultaneously. By Lemma \ref{LemHom}, this polynomial also annihilates $f$, which yields the desired decomposition.
\end{proof}

\begin{theorem}[Krull-Schmidt]
For $M \in \mathcal{M}_{\mathfrak{p}}$, there exists a finite decomposition of $M$ into indecomposable objects whose endomorphism rings are local. In other words, $\mathcal{M}_{\mathfrak{p}}$ is a Krull–Schmidt category. Consequently, this decomposition into indecomposables is unique up to permutation and isomorphism of the factors.
\end{theorem}
\newpage
\begin{proof} \mbox{}
\begin{itemize}
\item If $M \in \mathcal{M}_{\mathfrak{p}}$ is not indecomposable, we can write $M = M_0 \oplus M_1$ and then recursively decompose each direct summand. This process must terminate after finitely many steps by Lemma \ref{LemDecomposition}, since every nonzero submodule contains at least one maximal vector. Hence the number of direct summands is bounded by $\sum_{\lambda \in \mathfrak{h}^*} \dim M_\lambda^+$, which is finite.
\item If $M = \bigoplus_{i=1}^n M_i$ is a decomposition into indecomposables, then for each $i$, the endomorphism ring $\operatorname{End}_{\mathfrak{b}}(M_i)$ is local by Lemma \ref{LemFitting}.
\end{itemize}
Therefore, $\mathcal{M}_{\mathfrak{p}}$ is a Krull–Schmidt category. Consequently, the indecomposable summands appearing in a decomposition of $M \in \mathcal{M}_{\mathfrak{p}}$ into indecomposables are unique up to permutation and isomorphism of the factors (see e.g. \cite[Thm.~4.2]{Krause2015}).
\end{proof}

\subsection{Filtration of direct summands}

For $M \in \mathcal{O}$, the module $M$ admits a finite filtration whose subquotients are simple modules $L(\mu)$. It turns out that an analogous type of filtration exists for the direct summands of $\operatorname{Res}_{\mathfrak{p}}^{\mathfrak{g}} M$.

\begin{proposition}
\label{PropFiltration}
Let $M \in \mathcal{O}$ and let $A$ be a direct summand of $\operatorname{Res}_{\mathfrak{p}}^{\mathfrak{g}} M$. Then $A$ admits a filtration:
$$
0 = A_0 \subset A_1 \subset \cdots \subset A_n = A
$$
such that $A_{i+1} / A_i \cong \operatorname{Res}_{\mathfrak{p}}^{\mathfrak{g}} L(\lambda_i)$ for some $\lambda_i \in \mathfrak{h}^*$. \\

More precisely, if $\operatorname{Res}_{\mathfrak{p}}^{\mathfrak{g}} M = \bigoplus_{j=1}^m A^{(j)}$, then there exists a composition series of $M$ and filtrations of the $A^{(j)}$ as above, compatible in the following sense: for each $i$ there exist indices $a_j$ and $k$ such that
\begin{align*}
& \operatorname{Res}_{\mathfrak{p}}^{\mathfrak{g}} M_i = A^{(1)}_{a_1} \oplus \cdots \oplus A^{(n)}_{a_n}, \\
& \operatorname{Res}_{\mathfrak{p}}^{\mathfrak{g}} M_{i+1} = A^{(1)}_{a_1} \oplus \cdots \oplus A^{(k)}_{a_k+1} \oplus \cdots \oplus A^{(n)}_{a_n},
\end{align*}
with $M_{i+1}/M_i \cong L(\mu)$ and $A^{(k)}_{a_k+1} / A^{(k)}_{a_k} \cong \operatorname{Res}_{\mathfrak{p}}^{\mathfrak{g}} L(\mu)$.
\end{proposition}
\begin{proof}
For simplicity, we give the argument for the case of two direct summands $\operatorname{Res}_{\mathfrak{p}}^{\mathfrak{g}} M = A \oplus B$. We prove the existence of a filtration of $A$ and $B$ by induction on the length of $M$. Let $v^+$ be a maximal vector of $M$ of minimal weight $\mu$ with respect to $\leq$. By Lemma \ref{LemDecomposition}, we may assume that $v^+ \in A$ or $v^+ \in B$. Without loss of generality, assume $v^+ \in A$. Since $v^+$ is a maximal vector of minimal weight, we have an inclusion $\iota : L(\mu) \hookrightarrow M$. The image of $\iota$ is entirely contained in $A$. Indeed, if $f(v) = a + b \in A \oplus B$ for $v \in L(\mu)$ with $b \neq 0$, then there exists $x \in U(\mathfrak{b})$ such that $x \cdot b$ is a maximal vector of $M$. By the minimality of $\mu$, we have $x \cdot v = 0$, hence $x \cdot a = - x \cdot b \neq 0$, a contradiction. We conclude that
$$
\operatorname{Res}_{\mathfrak{p}}^{\mathfrak{g}} M / L(\mu) = A / L(\mu) \oplus B.
$$
By the induction hypothesis, $A / L(\mu)$ and $B$ admit filtrations of the desired form, which therefore also yields a filtration for $A$.
\end{proof}

\subsection{Split morphisms and decomposition results}

The functors $(\operatorname{Ind}^\mathfrak{g}_\mathfrak{p}, \operatorname{Res}^\mathfrak{g}_\mathfrak{p})$ form an adjoint pair between the category of $\mathfrak{g}$-modules and the category of $\mathfrak{p}$-modules. We will use an elementary result concerning adjoint functors and split morphisms. Let $\mathcal{C}, \mathcal{D}$ be two categories, and let $A, B$ be objects in $\mathcal{C}$ or $\mathcal{D}$.

\begin{definition}
If $f: A \to B$ and $g: B \to A$ satisfy $f \circ g = \operatorname{id}_B$, then $f$ is called a split epimorphism and $g$ a split monomorphism.
\end{definition}

\begin{lemma}
\label{LemSplit}
Let $(L, R)$ be an adjoint pair of functors with $L: \mathcal{C} \to \mathcal{D}$ and $R: \mathcal{D} \to \mathcal{C}$. Then the following equivalences hold:
\begin{enumerate}
\item $LA \xrightarrow{Lf} LB$ is split mono if and only if $RLA \xrightarrow{RLf} RLB$ is split mono.
\item $LA \xrightarrow{f} LB$ is split epi if and only if $RLA \xrightarrow{Rf} RLB$ is split epi.
\item $RA \xrightarrow{f} RB$ is split mono if and only if $LRA \xrightarrow{Lf} LRB$ is split mono.
\item $RA \xrightarrow{Rf} RB$ is split epi if and only if $LRA \xrightarrow{LRf} LRB$ is split epi.
\end{enumerate}
\end{lemma}

\begin{proof}
Since all functors send split morphisms to split morphisms, the forward implications are immediate. Let $\eta: \operatorname{Id} \to RL$ be the unit and $\varepsilon: LR \to \operatorname{Id}$ the counit of the adjunction.

\begin{itemize}
\item Suppose $RLA \xrightarrow{RLf} RLB$ is split mono. Then so is $LRLA \xrightarrow{LRLf} LRLB$. Let $r$ be a retraction of this morphism. Using the properties of adjoint functors, we obtain the following commutative diagram:

\begin{center}
$\begin{tikzcd}[sep=3em]
LA \arrow[r, "Lf"] \arrow[d, "L\eta"] \arrow[dd, bend right=50, "id"'] & LB \arrow[d, "L\eta"'] \\ 
LRLA \arrow[d, "\varepsilon L"] \arrow[r, hook, "LRLf"'] & LRLB \arrow[l, bend right=20, "r"'] \\
LA & 
\end{tikzcd}$
\end{center}
Then $\tilde{r} := (\varepsilon L) \circ r \circ (L\eta)$ is a retraction of $LA \xrightarrow{Lf} LB$.

\item Suppose $RLA \xrightarrow{Rf} RLB$ is split epi. Then so is $LRLA \xrightarrow{LRLf} LRLB$. Let $s$ be a section of this morphism. We then have the commutative diagram:

\begin{center}
$\begin{tikzcd}[sep=3em]
& LB \arrow[d, "L\eta"'] \arrow[dd, bend left=50, "id"] \\
LRLA \arrow[d, "\varepsilon L"] \arrow[r, two heads, "LRf"'] & LRLB \arrow[d, "\varepsilon L"'] \arrow[l, bend right=20, "s"'] \\
LA \arrow[r, "f"] & LB
\end{tikzcd}$
\end{center}

Then $\tilde{s} := (\varepsilon L) \circ s \circ (L\eta)$ is a section of $LA \xrightarrow{f} LB$.
\end{itemize}

The remaining two implications follow by duality.
\end{proof}
\noindent
As a consequence of Lemma \ref{LemSplit}, we obtain two types of direct summands of $\operatorname{Res}_{\mathfrak{p}}^{\mathfrak{g}} M$, for $M \in \mathcal{O}^\mathfrak{p}$: $\mathfrak{p}$-tilting submodules and finite-dimensional quotients.

\begin{proposition}
\label{PropDecompTilting}
Let $M \in \mathcal{O}^\mathfrak{p}$ and suppose $M$ contains a $\mathfrak{p}$-tilting submodule $N$. Then $\operatorname{Res}_{\mathfrak{p}}^{\mathfrak{g}} N$ is a direct summand of $\operatorname{Res}_{\mathfrak{p}}^{\mathfrak{g}} M$.
\end{proposition}

\begin{proof}
By exactness of the functor $\left(M_{\mathfrak{p}}(0) \otimes \cdot\right)^\chi$, the short exact sequence:
\begin{center}
$\begin{tikzcd}
N \arrow[hook, r] 
& M \arrow[two heads, r]
& M/N
\end{tikzcd}$
\end{center}
induces, for any central character $\chi$, a short exact sequence in $\mathcal{O}^\mathfrak{p}$:
\begin{center}
$\begin{tikzcd}
\left(M_{\mathfrak{p}}(0) \otimes N\right)^\chi \arrow[hook, r] 
& \left(M_{\mathfrak{p}}(0) \otimes M\right)^\chi \arrow[two heads, r]
& \left(M_{\mathfrak{p}}(0) \otimes \left( M/N \right)\right)^\chi
\end{tikzcd}$
\end{center}
By Proposition \ref{PropPreliFiltration}, the left term is a $\mathfrak{p}$-tilting module, while the right term has a generalized Verma filtration. Properties of tilting modules imply that this exact sequence splits for all $\chi$. Hence, $M_{\mathfrak{p}}(0) \otimes N$ is a direct summand of $M_{\mathfrak{p}}(0) \otimes M$, i.e. there is a split monomorphism
$\operatorname{Ind}_{\mathfrak{p}}^{\mathfrak{g}} \operatorname{Res}_{\mathfrak{p}}^{\mathfrak{g}} N \hookrightarrow \operatorname{Ind}_{\mathfrak{p}}^{\mathfrak{g}} \operatorname{Res}_{\mathfrak{p}}^{\mathfrak{g}} M$. By Lemma \ref{LemSplit}, this implies a split monomorphism $\operatorname{Res}_{\mathfrak{p}}^{\mathfrak{g}} N \hookrightarrow \operatorname{Res}_{\mathfrak{p}}^{\mathfrak{g}} M$, so $\operatorname{Res}_{\mathfrak{p}}^{\mathfrak{g}} N$ is indeed a direct summand of $\operatorname{Res}_{\mathfrak{p}}^{\mathfrak{g}} M$.
\end{proof}

\begin{corollary}
\label{CorDecompTilting}
Let $\lambda, \mu \in \Lambda^+_I$. If $M_{\mathfrak{p}}(\mu)$ is a $\mathfrak{p}$-tilting submodule of $M_{\mathfrak{p}}(\lambda)$, then as a $\mathfrak{p}$-module:
$$
M_{\mathfrak{p}}(\lambda) \cong M_{\mathfrak{p}}(\mu) \oplus \left(M_{\mathfrak{p}}(\lambda)/M_{\mathfrak{p}}(\mu)\right)
$$
where $M_{\mathfrak{p}}(\mu) = L(\mu)$ is indecomposable as a $\mathfrak{p}$-module. In particular, for any $\lambda' \in \Lambda^+_I$:
$$
M_{\mathfrak{p}}(\lambda) \otimes M_{\mathfrak{p}}(\lambda') \cong 
M_{\mathfrak{p}}(\mu) \otimes M_{\mathfrak{p}}(\lambda') \;\oplus\; 
\left(M_{\mathfrak{p}}(\lambda)/M_{\mathfrak{p}}(\mu)\right) \otimes M_{\mathfrak{p}}(\lambda')
$$
where $M_{\mathfrak{p}}(\mu) \otimes M_{\mathfrak{p}}(\lambda')$ decomposes into direct sums of indecomposable $\mathfrak{p}$-tilting modules $D_{\mathfrak{p}}(\nu)$.
\end{corollary}

\begin{remark}
The second part of the corollary can also be obtained directly by noting, via the same arguments as in Proposition \ref{PropDecompTilting}, that the following short exact sequence splits:
\begin{center}
$\begin{tikzcd}
M_{\mathfrak{p}}(\mu) \otimes M_{\mathfrak{p}}(\lambda') \arrow[hook, r] 
& M_{\mathfrak{p}}(\lambda) \otimes M_{\mathfrak{p}}(\lambda') \arrow[two heads, r]
& \left(M_{\mathfrak{p}}(\lambda)/M_{\mathfrak{p}}(\mu) \right)\otimes M_{\mathfrak{p}}(\lambda')
\end{tikzcd}$
\end{center}
\end{remark}

\begin{proposition}
If $M \in \mathcal{O}^\mathfrak{p}$ has a finite-dimensional quotient $V$, then $\operatorname{Res}_{\mathfrak{p}}^{\mathfrak{g}} V$ is a direct summand of $\operatorname{Res}_{\mathfrak{p}}^{\mathfrak{g}} M$.
\end{proposition}
\begin{proof}
Consider the short exact sequence:
\begin{center}
$\begin{tikzcd}
M_{\mathfrak{p}}(0) \otimes N \arrow[hook, r] 
& M_{\mathfrak{p}}(0) \otimes M \arrow[two heads, r]
& M_{\mathfrak{p}}(0) \otimes V
\end{tikzcd}$
\end{center}
It splits because $M_{\mathfrak{p}}(0)$ is projective in $\mathcal{O}^\mathfrak{p}$ and $V$ is finite-dimensional, hence $M_{\mathfrak{p}}(0) \otimes V$ is projective in $\mathcal{O}^\mathfrak{p}$. Therefore, $M_{\mathfrak{p}}(0) \otimes V$ is a direct summand of $M_{\mathfrak{p}}(0) \otimes M$, i.e. there is a split epimorphism
$\operatorname{Ind}_{\mathfrak{p}}^{\mathfrak{g}} \operatorname{Res}_{\mathfrak{p}}^{\mathfrak{g}} M \twoheadrightarrow \operatorname{Ind}_{\mathfrak{p}}^{\mathfrak{g}} \operatorname{Res}_{\mathfrak{p}}^{\mathfrak{g}} V$. By Lemma \ref{LemSplit}, this induces a split epimorphism $\operatorname{Res}_{\mathfrak{p}}^{\mathfrak{g}} M \twoheadrightarrow \operatorname{Res}_{\mathfrak{p}}^{\mathfrak{g}} V$.
\end{proof}

\begin{corollary}
\label{CorDecompHead}
For $\lambda \in \Lambda^+$, as a $\mathfrak{p}$-module:
$$
M_{\mathfrak{p}}(\lambda) \cong N_{\mathfrak{p}}(\lambda) \oplus L(\lambda)
$$
where $L(\lambda)$ is indecomposable as a $\mathfrak{p}$-module. In particular, for any $\lambda' \in \Lambda^+_I$:
$$
M_{\mathfrak{p}}(\lambda) \otimes M_{\mathfrak{p}}(\lambda') \cong 
N_{\mathfrak{p}}(\lambda) \otimes M_{\mathfrak{p}}(\lambda') \;\oplus\; 
L(\lambda) \otimes M_{\mathfrak{p}}(\lambda')$$
\end{corollary}

The following proposition provides a more explicit description of the complement of $N_{\mathfrak{p}}(\lambda)$ in $M_{\mathfrak{p}}(\lambda)$. In particular, it is unique.

\begin{proposition}
\label{PropDecompHead}
Let $\lambda \in \Lambda^+$ and $\mathfrak{p} = \mathfrak{p}_I$. Set $n_i = \langle -w_0 \lambda + \rho, \alpha_i^\vee \rangle$ for $i \in \{1, \dots, l\}$. Then $N_{\mathfrak{p}}(\lambda)$ admits a unique complement $L$ in $M_{\mathfrak{p}}(\lambda)$. This complement is generated, as a $\mathfrak{p}$-module, by nonzero vectors $v$ such that:
\begin{itemize}
\item $v \in M_{\mathfrak{p}}(\lambda)_{w_0 \lambda}$,
\item $x_i^{n_i} v = 0$ for all $i \in \{1, \dots, l\}$,
\item $y_j v = 0$ for all $j \in I$
\end{itemize}
In particular, such a vector $v$ is unique up to scalar multiplication. 
\end{proposition}
\begin{proof}
We have $\dim L(\lambda)_{w_0\lambda} = \dim L(\lambda)_\lambda = 1$. We first determine the annihilator of a nonzero $v \in L(\lambda)_{w_0 \lambda}$.
\begin{itemize}
\item We have $L(-w_0 \lambda) \cong U(\mathfrak{g}) / I_0$ where $$I_0  =\langle h - (-w_0 \lambda)(h) \cdot 1, \, y_i^{n_i}, \, \mathfrak{n} \mid  h \in \mathfrak{h},\, i \in \{1, \dots, l\} \rangle $$
(see \cite[§2.6]{Humphreys2008}).
\item Hence $L(\lambda) \cong L(-w_0 \lambda)^* \cong U(\mathfrak{g}) / I$ where 
$$I = \langle h - (-w_0 \lambda)(h) \cdot 1, \, x_i^{n_i}, \, \mathfrak{n}^- \mid h \in \mathfrak{h}, \, i \in \{1, \dots, l\}\rangle $$
using the fact that $L(-w_0 \lambda)^\vee \cong L(-w_0 \lambda)$ and the $\mathfrak{g}$-module structure on the dual is twisted by an anti-involution $\tau$ that fixes $\mathfrak{h}^*$ pointwise and send $\mathfrak{g}_\alpha$ to $\mathfrak{g}_{-\alpha}$.
\item Finally $\operatorname{Res}^{\mathfrak{g}}_\mathfrak{p} L(\lambda) \cong U(\mathfrak{p}) / J$ where 
$$J = \langle h - (-w_0 \lambda)(h) \cdot 1, \, x_i^{n_i}, \, y_j \mid h \in \mathfrak{h}, \, i \in \{1, \dots, l\}, \, j \in I \rangle$$ by the PBW theorem.
\end{itemize}

\begin{enumerate}
\item Let $L$ be a complement of $N_{\mathfrak{p}}(\lambda)$. Then $L \cong \operatorname{Res}^{\mathfrak{g}}_\mathfrak{p} L(\lambda)$, and the annihilator of a nonzero $v \in L_{w_0 \lambda}$ satisfies the three points of the proposition and generates $L$ as a $\mathfrak{p}$-module.

\item Conversely, if $v$ satisfies the conditions of the proposition, then $v \in M(\lambda)_{w_0 \lambda}$ implies $Ann_\mathfrak{p}(v) \subset J$, so $L(\lambda) = U(\mathfrak{p}) / J$ surjects onto $L = U(\mathfrak{p}) \cdot v = U(\mathfrak{b}) / Ann_\mathfrak{p}(v)$. Then for $w \in W$ with $w \neq \operatorname{id}$:
\begin{align*}
\operatorname{dim} L_{w\cdot\lambda} \leq \operatorname{dim} L(\lambda)_{(\lambda+\rho) - \rho} = \operatorname{dim} L(\lambda)_{\lambda+\rho -w(\rho)} = 0
\end{align*}
because the character of $L(\lambda)$ vanishes on weights $\mu \not\leq \lambda$, and $\rho \leq w(\rho)$ only if $w = \operatorname{id}$. Hence $L$ contains no maximal vector of $N_{\mathfrak{p}}(\lambda)$. By Lemma \ref{LemDecomposition}, $L \oplus N_{\mathfrak{p}}(\lambda)$ is direct. Let $\pi: M_{\mathfrak{p}}(\lambda) \to L(\lambda)$ be the projection with kernel $N_{\mathfrak{p}}(\lambda)$. Then $\pi: L \hookrightarrow L(\lambda)$ is injective, and since $L(\lambda)$ is generated by any of its vectors of weight $w_0 \lambda$, we have $L \cong L(\lambda)$, giving $M_{\mathfrak{p}}(\lambda) = N_{\mathfrak{p}}(\lambda) \oplus L$.

\item To prove uniqueness, note that $N_{\mathfrak{p}}(\lambda)_{w_0 \lambda}$ has codimension 1 in $M_{\mathfrak{p}}(\lambda)_{w_0 \lambda}$, and any vector $v \in M(\lambda)_{w_0 \lambda}$ in the intersection of the kernels of $x_i^{n_i}$ and $y_j$ is complementary to $N_{\mathfrak{p}}(\lambda)_{w_0 \lambda}$. This intersection is 1-dimensional, so the complement is unique. \qedhere
\end{enumerate}
\end{proof}

\section{Decomposition in the cases $\mathfrak{sl}_2$ and $\mathfrak{sl}_3$}
\label{Section4}

In this section, we provide the decomposition into indecomposable $\mathfrak{b}$-modules for certain modules in the category $\mathcal{O}$ in the cases of $\mathfrak{sl}_2$ and $\mathfrak{sl}_3$.

\subsection{Case $\mathfrak{g} = \mathfrak{sl}_2$}

Let $\mathfrak{g} = \mathfrak{sl}_2$, and identify $\mathfrak{h}^*$ with $\mathbb{C}$ so that the unique root corresponds to $2$. We fix $\lambda \geq -1$ and set $\mu = -\lambda -2$. The indecomposable modules in $\mathcal{O}_{\chi_\lambda}$ are (see \cite[Prop.~3.12]{Humphreys2008}):
$$
L(\lambda), \quad M(\mu)= L(\mu)=D(\mu), \quad M(\lambda)=P(\lambda), \quad M(\lambda)^\vee, \quad P(\mu)=D(\lambda).
$$
All of these coincide when $\lambda \notin \mathbb{N}$ and are pairwise non-isomorphic if $\lambda \in \mathbb{N}$. From now on, we consider this latter case.

\begin{proposition}
We have the following decompositions into indecomposable $\mathfrak{b}$-modules:
\begin{itemize}
\item $M(\lambda) = L(\lambda) \oplus M(\mu)$,
\item $P(\mu) = L(\lambda) \oplus 2\, M(\mu)$.
\end{itemize}
The $\mathfrak{b}$-modules $L(\lambda)$, $M(\mu)$ and $M(\lambda)^\vee$ are indecomposable.
\end{proposition}
\begin{proof}\mbox{}
\begin{itemize}
\item By Corollary \ref{CorDecompTilting} or \ref{CorDecompHead}, the decomposition of $M(\lambda)$ holds. 
\item For $P(\mu)$, consider the short exact sequence of $\mathfrak{b}$-modules:
$$
\begin{tikzcd}
M(\mu) \arrow[hook, r] & P(\mu) \arrow[two heads, r] & M(\mu)\oplus L(\lambda)
\end{tikzcd}
$$
This sequence splits by Corollary \ref{CorDecompTilting}, yielding the decomposition of $P(\mu)$.
\item It remains to show that $M(\lambda)^\vee$ is indecomposable. Suppose, by contradiction, that $M(\lambda)^\vee$ is decomposable. Then, by Proposition \ref{PropFiltration}, we would have the decomposition
$M(\lambda)^\vee = L(\lambda) \oplus M(\mu)$ as $\mathfrak{b}$-modules. However, $L(\lambda) \subset M(\lambda)^\vee$ is a $\mathfrak{g}$-submodule, and $M(\mu) \subset M(\lambda)^\vee$ is stable by $\mathfrak{n}^-$ by considering the weight spaces. This would yield a decomposition of $\mathfrak{g}$-module $M(\lambda)^\vee = L(\lambda) \oplus M(\mu)$ which is a contradiction. \qedhere
\end{itemize}
\end{proof}

Thanks to the decomposition into $\mathfrak{b}$-modules, we can easily obtain the decomposition into indecomposables of all tensor products $M(\lambda') \otimes M$ in $\mathfrak{sl}_2$ by noting that:
\begin{itemize}
\item $M(\lambda') \otimes L(\lambda)$ is projective for dominant $\lambda'$ and $\lambda \in \mathbb{N}$, hence decomposes as a direct sum of indecomposable projectives $M(\lambda'')$ and $P(\mu'')$,
\item $M(\lambda') \otimes L(\lambda)$ is tilting for antidominant $\lambda'$ and $\lambda \in \mathbb{N}$, hence decomposes as a direct sum of indecomposable tiltings $M(\mu'')$ and $P(\mu'')$,
\item $M(\lambda') \otimes M(\lambda)$ decomposes into a sum of indecomposable tiltings if either $\lambda$ or $\lambda'$ is antidominant, by Proposition \ref{PropPreliFiltration},
\item $M(\lambda') \otimes M(\lambda)^\vee$ decomposes into a sum of indecomposable tiltings, again by Proposition \ref{PropPreliFiltration}.
\end{itemize}
Since the decomposition of projectives and tiltings is entirely determined by their character, it can be computed explicitly. For instance, we have:
\begin{align*}
M(0)\otimes M(0) &=  M(0) \otimes L(0) \,\oplus \,  M(0)\otimes  M(-2) \\
&= M(0) \oplus \bigoplus_{i \geq 1} M(-2i) \\
M(0) \otimes M(0)^\vee &= P(-2) \oplus \bigoplus_{i \geq 2} M(-2i) \\
M(0)\otimes P(-2) &=  M(0)\otimes L(0) \,\oplus \,   2M(0)\otimes M(-2) \\
&= M(0) \oplus \bigoplus_{i \geq 1} 2M(-2i) 
\end{align*}
Note that the problem of explicitly computing the decomposition of the tensor product of two Verma modules was already solved in the case of $\mathfrak{sl}_2$ in \cite{Murakami2024}.

\subsection{Case $\mathfrak{g}=\mathfrak{sl}_3$}

Let $\mathfrak{g} = \mathfrak{sl}_3$ with a chosen Chevalley basis:
$$
y_\alpha, y_\beta, y_{\alpha+\beta} := [y_\alpha, y_\beta], \quad h_\alpha, h_\beta, \quad x_\alpha, x_\beta, x_{\alpha+\beta} := [x_\alpha, x_\beta]
$$
We will now provide the decomposition of the Verma modules in $\mathcal{O}_{\chi_\lambda}$ into $\mathfrak{b}$-modules, for $\lambda \in \Lambda^{+}$. Let $\lambda \in \Lambda^{+}$, and set: 
$$n = \langle \lambda+\rho, \alpha^\vee \rangle \in \mathbb{Z}^{\geq 1},\quad m = \langle \lambda+\rho, \beta^\vee \rangle\in \mathbb{Z}^{\geq 1}$$ 
The Verma module $M(\lambda)$ has 5 proper Verma submodules corresponding to the 5 other linearly independant maximal vectors: 
\begin{itemize}
\item $M(s_\alpha\cdot \lambda) = U(\mathfrak{n}^-)v_\alpha$ where $v_\alpha=y_\alpha^{n}\cdot v^{+}$,
\item $M(s_\beta\cdot \lambda) = U(\mathfrak{n}^-)v_\beta$ where $v_\beta=y_\beta^{m} \cdot v^{+}$,
\item $M(s_\beta s_\alpha\cdot \lambda) = U(\mathfrak{n}^-)v_{\beta\alpha}$ where $v_{\beta\alpha}=y_\beta^{n+m}y_\alpha^{n}\cdot v^{+}$,
\item $M(s_\alpha s_\beta\cdot \lambda) = U(\mathfrak{n}^-)v_{\alpha\beta}$ where $v_{\alpha\beta}=y_\alpha^{n+m}y_\beta^{m}\cdot v^{+}$,
\item $M(w_0\cdot \lambda) = U(\mathfrak{n}^-)v_0$ where $v_0=y_\alpha^{m}y_\beta^{n+m}y_\alpha^{n}.v^{+}\propto y_\beta^{n}y_\alpha^{n+m}y_\beta^{m}\cdot v^{+}$
\end{itemize}
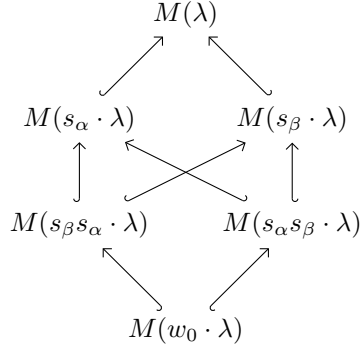
\begin{center}
\begin{tikzpicture}[scale=1.4, hook/.style={{Hooks[right]}-{Straight Barb}},  hook'/.style={{Hooks[left]}-{Straight Barb}}]
	\node (e) at (0,3)    {$M(\lambda)$};
	\node (s1) at (-1, 2) {$M(s_\alpha\cdot\lambda)$};
	\node (s2) at (1, 2) {$M(s_\beta\cdot\lambda)$};
	\node (s1s2) at (-1, 1) {$M(s_\beta s_\alpha\cdot\lambda)$};
	\node (s2s1) at (1, 1) {$M(s_\alpha s_\beta \cdot\lambda)$};
	\node (w0) at (0, 0) {$M(w_0\cdot \lambda)$};
\draw (s1) edge[hook] node[] {} (e);
\draw (s2) edge[hook'] node[] {} (e);
\draw (s1s2)  edge[hook] node[] {} (s1);
\draw  (s2s1) edge[hook'] node[] {} (s2);
\draw  (s2s1) edge[hook'] node[] {} (s1);
\draw  (s1s2) edge[hook] node[] {} (s2);
\draw (w0)  edge[hook] node[] {} (s1s2);
\draw (w0)  edge[hook' ] node[] {} (s2s1);
\end{tikzpicture}
\end{center}

\begin{proposition}
\label{PropSL3}
For $\mathfrak{g}= \mathfrak{sl}_3$ and $\lambda \in \Lambda^{+}$, the Verma module $M(\lambda)$ admits the decomposition into indecomposable $\mathfrak{b}$-modules as follows: 
\begin{align*}
M(\lambda) = L(\lambda) \oplus \bigl( M(s_{\alpha}\cdot \lambda)/ M(s_{\alpha}s_\beta \cdot \lambda) \bigr) \oplus \bigl( M(s_{\beta}\cdot \lambda)/ M(s_\beta s_\alpha \cdot \lambda) \bigr) \oplus M(w_0 \cdot \lambda)
\end{align*}
\end{proposition}

\begin{proof}
According to Proposition \ref{PropDecompHead}, we have the decomposition $M(\lambda) = L(\lambda) \oplus N(\lambda)$. The strategy is to decompose $N(\lambda)$ by constructing suitable submodules and applying Lemma \ref{LemDecomposition}. We will then establish the announced $\mathfrak{b}$-module isomorphisms and verify that these modules are indeed indecomposable. Consider the following three $\mathfrak{b}$-submodules of $N(\lambda)$:

\begin{enumerate}
\item $M_\alpha \subset M(s_\alpha \cdot \lambda)$, the vector space generated by $y_{\alpha+\beta}^l y_\beta^p y_\alpha^k y_\alpha^n \cdot v^+$ for $p,k \ge 0$ and $l < m$,
\item $M_\beta \subset M(s_\beta \cdot \lambda)$, generated by $y_{\alpha+\beta}^l y_\alpha^p y_\beta^k y_\beta^m \cdot v^+$ for $p,k \ge 0$ and $l < n$,
\item $M(w_0 \cdot \lambda)$
\end{enumerate}

\begin{itemize}
\item Claim 1) $M_\alpha, M_\beta$ are $\mathfrak{b}$-modules.  Since $y_{\alpha+\beta}$ commutes with both $y_\alpha$ and $y_\beta$, we may position the $y_{\alpha+\beta}$ factors between the $y_\alpha$ and $y_\beta$ in the definitions of $M_\alpha$ and $M_\beta$. In this form, both subspaces are stable under the actions of $x_\alpha$ and $x_\beta$.

\item Claim 2) $M_\alpha \oplus M_\beta \oplus M(w_0 \cdot \lambda) = N(\lambda)$. We first verify that the vectors $v_\alpha, v_{\beta\alpha} \in M_\alpha$ are not in $M_\beta$, which is immediate for $v_\alpha$. For $v_{\beta\alpha}$, an inductive calculation shows that the commutator $[y_\alpha^i, y_\beta^j]$ contains a term  $y_{\alpha+\beta}^i y_\beta^{j-i}$ for $i \ge j$. In particular, $[y_\alpha^n, y_\beta^{n+m-1}]$ produces $y_{\alpha+\beta}^n y_\beta^{m-1}$, which, thanks to the condition $l < n$, ensures $v_{\beta\alpha} \notin M_\beta$. By symmetry, the vectors $v_\beta, v_{\alpha\beta} \in M_\beta$ are not in $M_\alpha$. Similarly, $v_0$ lies outside both $M_\alpha$ and $M_\beta$. As the character of $M(w_0 \cdot \lambda)$ vanishes on the weights of $v_\alpha, v_\beta, v_{\alpha\beta}, v_{\beta\alpha}$, the maximal vectors of $N(\lambda)$ are contained in exactly one of the three submodules. Lemma \ref{LemDecomposition} then yields
$M_\alpha \oplus M_\beta \oplus M(w_0 \cdot \lambda)$, which is equal to $N(\lambda)$ by weight-space dimensions.

\item Claim 3) $M_\alpha \cong M(s_\alpha \cdot \lambda) / M(s_\alpha s_\beta \cdot \lambda)$. According to the proof of Proposition \ref{PropFiltration}, we deduce from the decomposition of $N(\lambda)$ the decomposition $M(s_\alpha \cdot \lambda) = M_\alpha \oplus L(s_\alpha s_\beta \cdot \lambda) \oplus M(w_0 \cdot \lambda)$. Quotienting first by $M(w_0 \cdot \lambda) = L(w_0 \cdot \lambda)$ and then by $L(s_\alpha s_\beta \cdot \lambda)$, we obtain $M(s_\alpha \cdot \lambda) / M(s_\alpha s_\beta \cdot \lambda) = M_\alpha$.

\item Claim 4) $M_\beta \cong M(s_\beta \cdot \lambda) / M(s_\beta s_\alpha \cdot \lambda)$.
By the symmetry $\alpha \leftrightarrow \beta$.

\item Claim 5) $M_\alpha, M_\beta$ indecomposable.
Let us show that $M_\alpha$ is indecomposable. Consider the map 
$$
x_\alpha^{m} : (M_\alpha)_{w_0 \lambda} \longrightarrow (M_\alpha)_{\lambda - (n+m)\beta - n\alpha}
$$
We claim that this map is nonzero, which ensures that one of the vectors in $(M_\alpha)_{w_0 \lambda}$ is sent to $v_{\beta\alpha}$. We also claim that the map is injective, which implies that all nonzero vectors in $(M_\alpha)_{w_0 \lambda}$ can be sent to $v_\alpha$. Hence, if $M_\alpha$ were decomposable into a module containing $v_\alpha$ and another containing $v_{\beta\alpha}$, the first module would necessarily contain the entire weight space $(M_\alpha)_{w_0 \lambda}$. But then $x_\alpha^m$ would send one of its vectors to $v_{\beta\alpha}$, which is a nonzero vector of the second module, a contradiction.

\item Claim 6) $x_\alpha^m$ nonzero. 
The map $x_\alpha^m$ sends $y_\alpha^{n+m} \cdot v^+$ to a nonzero multiple of $y_\alpha^n \cdot v^+$. Otherwise, one of the vectors $y_\alpha^{n+i} \cdot v^+$ would be maximal. Since $x_\alpha$ commutes with $y_\beta$, it follows that $x_\alpha^m$ sends $y_\beta^{n+m} y_\alpha^{n+m} \cdot v^+$ to a nonzero multiple of $y_\beta^{n+m} y_\alpha^n \cdot v^+ = v_{\beta\alpha}$.

\item Claim 7) $x_\beta$ injective. 
We compute the image of vectors in $(M_\alpha)_{w_0 \lambda}$ under $x_\beta$. For $i < m$:
\begin{align*}
& x_\beta (y_\beta^{n+m-i} y_{\alpha+\beta}^i y_\alpha^{n+m-i} \cdot v^{+}) \\
= \; & y_\beta^{n+m-i} x_\beta y_{\alpha+\beta}^i y_\alpha^{n+m-i} \cdot v^{+} + [x_\beta,y_\beta^{n+m-i}] y_{\alpha+\beta}^i y_\alpha^{n+m-i} \cdot v^{+}
\end{align*}
The first term is proportional to 
$$
y_\beta^{n+m-i} y_{\alpha+\beta}^{i-1} y_\alpha^{n+m-i+1} \cdot v^+
$$
if $i > 0$, and $0$ otherwise. The second term is
$$
C \, y_\beta^{n+m-i-1} y_{\alpha+\beta}^i y_\alpha^{n+m-i} \cdot v^+
$$
with $C = i - m < 0$. In particular, $C \neq 0$, showing that the image of $x_\beta$ has dimension $m$, hence $x_\beta$ is injective. This establishes the indecomposability of $M_\alpha$. The indecomposability of $M_\beta$ follows by the symmetry $\alpha \leftrightarrow \beta$.\qedhere
\end{itemize}
\end{proof}

The proof of the theorem also provides the decomposition of the other Verma modules in the orbit of $\lambda$ as $\mathfrak{b}$-modules:

\begin{proposition}
For $\mathfrak{g} = \mathfrak{sl}_3$ and $\lambda \in \Lambda^+$, we have the following decomposition into indecomposable $\mathfrak{b}$-modules:
\begin{align*}
M(s_\alpha \cdot \lambda) &= \bigl( M(s_\alpha \cdot \lambda)/M(s_\alpha s_\beta \cdot \lambda) \bigr) \oplus L(s_\alpha s_\beta \cdot \lambda) \oplus M(w_0 \cdot \lambda),\\
M(s_\beta \cdot \lambda) &= \bigl( M(s_\beta \cdot \lambda)/M(s_\beta s_\alpha \cdot \lambda) \bigr) \oplus L(s_\beta s_\alpha \cdot \lambda) \oplus M(w_0 \cdot \lambda),\\
M(s_\alpha s_\beta \cdot \lambda) &= L(s_\alpha s_\beta \cdot \lambda) \oplus M(w_0 \cdot \lambda),\\
M(s_\beta s_\alpha \cdot \lambda) &= L(s_\beta s_\alpha \cdot \lambda) \oplus M(w_0 \cdot \lambda)
\end{align*}
\end{proposition}

\begin{proof}
Corollary \ref{CorDecompTilting} and Proposition \ref{PropFiltration} yields the decomposition of $M(s_\alpha s_\beta \cdot \lambda)$.  In the proof of Proposition \ref{PropSL3}, we also established that $M(s_\alpha \cdot \lambda) = M_\alpha \oplus M(s_\beta s_\alpha \cdot \lambda)$, from which the first decomposition follows. The remaining decompositions are obtained by symmetry.
\end{proof}

Using the $\mathfrak{b}$-module decomposition, we can compute the decomposition into indecomposables of the tensor product of Verma modules in the principal block, restricted to the principal block. We adopt the notation $\{A,B\}$ to denote a module given by a non-split extension of $A$ by $B$.

\begin{enumerate}
\item $\begin{aligned}[t]
\bigl(M(0)\otimes M(0)\bigr)^{\chi_0} = & \; M(0) \\
\oplus & \; \bigl\lbrace M(w_0\cdot 0), M(s_\alpha s_\beta \cdot 0) \oplus M(s_\beta s_\alpha \cdot 0) \oplus M(s_\alpha \cdot 0) \bigr\rbrace \\
\oplus & \; \bigl\lbrace M(w_0\cdot 0), M(s_\alpha s_\beta \cdot 0) \oplus M(s_\beta s_\alpha \cdot 0) \oplus M(s_\beta \cdot 0) \bigr\rbrace \\
\oplus & \; M(w_0\cdot 0),
\end{aligned}$

\item $\begin{aligned}[t]
\bigl(M(0)\otimes M(s_\alpha \cdot 0)\bigr)^{\chi_0} = & \; M(s_\alpha \cdot 0) \\ \oplus & \; \bigl\lbrace M(w_0\cdot 0), M(s_\alpha s_\beta \cdot 0) \oplus M(s_\beta s_\alpha \cdot 0) \bigr\rbrace \\ 
\oplus & \; M(s_\alpha s_\beta \cdot 0) \\
\oplus & \; M(w_0\cdot 0),
\end{aligned}$

\item $\begin{aligned}[t] 
\bigl(M(0)\otimes M(s_\alpha s_\beta \cdot 0)\bigr)^{\chi_0} = & \; M(s_\alpha s_\beta \cdot 0) \\
\oplus & \; M(w_0\cdot 0),
\end{aligned}$

\item $\begin{aligned}[t] 
\bigl(M(0)\otimes M(w_0\cdot 0)\bigr)^{\chi_0} = & \; M(w_0\cdot 0),
\end{aligned}$

\item $\begin{aligned}[t] 
\bigl(M(s_\alpha \cdot 0)\otimes M(s_\alpha \cdot 0)\bigr)^{\chi_0} = & \; 
M(s_\alpha s_\beta \cdot 0) \\
\oplus & \; M(w_0\cdot 0),
\end{aligned}$

\item $\begin{aligned}[t] 
\bigl(M(s_\alpha \cdot 0)\otimes M(s_\beta \cdot 0)\bigr)^{\chi_0} = & \; \bigl\lbrace M(w_0\cdot 0), M(s_\alpha s_\beta \cdot 0)\oplus M(s_\beta s_\alpha \cdot 0) \bigr\rbrace \\
\oplus & \; M(w_0\cdot 0),
\end{aligned}$

\item $\begin{aligned}[t]
\bigl(M(s_\alpha \cdot 0)\otimes M(s_\beta s_\alpha \cdot 0)\bigr)^{\chi_0} = & \; M(w_0\cdot 0)
\end{aligned}$
\end{enumerate}
\noindent
The remaining tensor products in the principal block, which cannot be obtained by the symmetry $\alpha \leftrightarrow \beta$ from the seven previous cases, vanish.

\begin{proof} We begin with the first case, which is the only difficult one. From the decomposition of $\operatorname{Res}^{\mathfrak{g}}_{\mathfrak{p}} M(0)$ we obtain:
$$
M(0)\otimes M(0) = M(0) \,\oplus \, \operatorname{Ind}_{\mathfrak{p}}^{\mathfrak{g}} M_\alpha \,\oplus \, \operatorname{Ind}_{\mathfrak{p}}^{\mathfrak{g}} M_\beta \,\oplus \, M(0)\otimes M(w_0 \cdot \lambda)
$$
We have $\bigl(M(0)\otimes M(w_0 \cdot \lambda)\bigr)^{\chi_0} = M(w_0 \cdot \lambda)$, and by symmetry it suffices to show that $N = \bigl(\operatorname{Ind}_{\mathfrak{p}}^{\mathfrak{g}} M_\alpha\bigr)^{\chi_0}$ has the claimed form and is indecomposable. Since all nonzero weight spaces of $M_\alpha$ are one-dimensional, the module $N$ admits a Verma filtration with subquotients $M(w_0 \cdot 0), \; M(s_\beta s_\alpha \cdot 0),\; M(s_\alpha s_\beta \cdot 0), \; M(s_\alpha \cdot 0)$ each occurring with multiplicity one. Because $M_\alpha$ contains the maximal vectors $v_{\beta\alpha}, v_\alpha$, the module $N$ contains $M(s_\beta s_\alpha \cdot 0) \oplus M(s_\alpha \cdot 0)$ as submodules. A direct computation then shows that the subspace of maximal vectors $N^+_{s_\alpha s_\beta \cdot 0}$ is two-dimensional, and hence contains both the maximal vector associated with $M(s_\alpha s_\beta \cdot 0) \subset M(s_\alpha \cdot 0)$ and another linearly independent maximal vector, yielding a second copy of $M(s_\alpha s_\beta \cdot 0)$. Thus we indeed obtain $M(s_\alpha s_\beta \cdot 0) \oplus M(s_\beta s_\alpha \cdot 0) \oplus M(s_\alpha \cdot 0) \subset N,
$ which is the claimed form. \\

To establish the indecomposability of $N$, we first show that $N' := N / M(s_\alpha \cdot 0)$ is indecomposable. The module $N'$ contains $M(s_\beta s_\alpha \cdot 0) \oplus M(s_\alpha s_\beta \cdot 0)$ as submodule. Since the weight space $(M_\alpha)_{w_0 \cdot 0}$ is mapped bijectively onto $(M_\alpha)_{s_\beta s_\alpha \cdot 0}$ by $x_\alpha$, the quotient $N' / M(s_\beta s_\alpha \cdot 0)$ contains only one linearly independent maximal vector of weight $w_0 \cdot 0$. By contradiction, if $N' / M(s_\beta s_\alpha \cdot 0)$ were decomposable, it would split as $M(s_\beta s_\alpha \cdot 0) \oplus M(w_0 \cdot 0)$, and hence would contain two distinct submodules isomorphic to $M(w_0 \cdot 0)$, which is impossible. Similarly, $N'/M(s_\alpha s_\beta \cdot 0)$ is indecomposable, and therefore $N'$ itself is indecomposable. Now suppose, for contradiction, that $N$ is decomposable. Then one would have $N \cong N' \oplus M(s_\alpha \cdot 0)$, so that $N'' := N/\bigl(M(s_\alpha s_\beta \cdot 0) \oplus M(s_\beta s_\alpha \cdot 0)\bigr) = M(w_0 \cdot 0) \oplus M(s_\alpha \cdot 0)$. In particular, $N''$ would contain two maximal vectors of weight $w_0 \cdot 0$. A direct computation, however, shows that the space $(N'')^+_{w_0 \cdot 0}$ is one-dimensional, yielding a contradiction. \\

For the remaining cases, the decomposition of the tensor product follows directly from the $\mathfrak{b}$-module decomposition of one of the two factors. The indecomposability of modules of the form
$\{ M(w_0 \cdot 0), \; M(s_\alpha s_\beta \cdot 0)\oplus M(s_\beta s_\alpha \cdot 0) \}$ is proved by the same argument as above for the indecomposability of $N'$.
\end{proof}

\printbibliography
\end{document}